\documentclass[12pt]{amsart}
\usepackage{amsmath}
\usepackage{latexsym}
\usepackage{amssymb}
\usepackage{amsfonts}

\usepackage{graphics}
\include{PDF}

\renewcommand{\proof}{\par\noindent{\it Proof.\ \ }}
\def\qed{\ifmmode\square\else\nolinebreak\hfill
$\Box$\fi\par\vskip12pt}

\def\l{\langle} \def\r{\rangle}

    \def\ZZ{\mathbb Z}
 \def\CC{{\mathcal C}}

\def\mod{{\sf mod~}} \def\val{{\sf val}}

\def\soc{{\sf soc}}

\def\K{{\sf K}}

\def\Ga{{\it \Gamma}} \def\Sig{{\it \Sigma}} \def\Del{{\it \Delta}}
\def\Ome{{\it \Omega}}

\def\Cay{{\sf Cay}} \def\Cos{{\sf Cos}} 
\def\Aut{{\sf Aut}}  \def\Out{{\sf Out}}

   \def\CD{{\rm CD}}
\def\HS{{\rm HS}}

\def\D{{\sf D}} 
\def\S{{\sf S}} 
\def\J{{\sf J}} \def\M{{\sf M}}

\def\a{\alpha}  \def\d{\delta}

\def\PSp{{\sf PSp}}

\def\A{{\sf A}}\def\Sym{{\sf Sym}}

\def\PSL{{\sf PSL}}\def\PGL{{\sf PGL}}
 \def\SL{{\sf SL}}

\def\ASL{{\sf ASL}}

 \def\PSU{{\sf PSU}}  
 \def\F{{\sf F}} \def\D{{\sf D}}

\def\Sz{{\sf Sz}}

\newtheorem{theorem}{Theorem}[section]%
\newtheorem{lemma}[theorem]{Lemma}%
\newtheorem{corollary}[theorem]{Corollary}%
\newtheorem{example}[theorem]{Example}%

\input amssym.def
\input amssym.tex

\begin{document}

\title[symmetric graphs of valency seven]
{Symmetric graphs of valency seven and their basic normal quotient graphs}

\thanks{1991 MR Subject Classification 20B15, 20B30, 05C25.}
\thanks{This work was partially supported by
National Natural Science Foundation of China (11461007, 11231008).}

\author[Pan]{Jiangmin Pan}
\address{J. M. Pan\\
School of Statistics and Mathematics\\
Yunnan University of Finance and Economics\\
Kunming \\
P. R. China}
\email{jmpan@ynu.edu.cn}

\author[Huang]{Junjie Huang}
\address{J. J. Huang\\
School of Statistics and Mathematics\\
Yunnan University of Finance and Economics\\
Kunming \\
P. R. China}

\author[Wang]{Chao Wang}
\address{C. Wang\\
School of Statistics and Mathematics\\
Yunnan University of Finance and Economics\\
Kunming \\
P. R. China}

\maketitle

\begin{abstract}
A graph $\Ga$ is {\it basic} if $\Aut\Ga$
has no normal subgroup $N\ne 1$
such that $\Ga$ is a normal cover of the normal quotient graph $\Ga_N$.
In this paper, we completely determine
the basic normal quotient graphs
of all connected $7$-valent symmetric graphs of order $2pq^n$
with $p<q$ odd primes,
which consist of an infinite family of
dihedrants of order $2p$ with
$p\equiv 1(\mod 7)$,
and 6 specific graphs with order at most $310$.
As a consequence,
it shows that, for any given positive integer $n$, there are only finitely many
connected $2$-arc-transitive $7$-valent graphs of order $2pq^n$ with $7\ne p<q$ primes,
partially generalizing Theorem 1 of
Conder, Li and Poto\v cnik [On the orders of arc-transitive graphs,
{\it J. Algebra} {\bf 421} (2015), 167--186].
\end{abstract}

\qquad {\textsc k}{\scriptsize \textsc {eywords.}} {\footnotesize
symmetric graph, basic normal quotient graph, normal cover}

\section{Introduction}

In this paper, graphs are undirected and
have no loops and multiple edges.
For a graph $\Ga$,
we denote by $V\Ga$ and $A\Ga$ its vertex set and arc set
respectively, and by $\Aut\Ga$ its full automorphism group.
The size $|V\Ga|$ is called the {\it order} of $\Ga$.
If there is a group $G\le\Aut\Ga$ acting transitively on $V\Ga$ or $A\Ga$, then
$\Ga$ is called  $G$-{\it vertex-transitive}
or $G$-{\it arc-transitive}, respectively.
An arc-transitive graph is also called a {\it symmetric graph}.

In the field of algebraic combinatorics,
symmetric graphs with order a fixed number (generally small) times a prime power have received
a lot of attention.
Chao \cite{Chao}, Cheng and Oxley \cite{C-O}
and Wang and Xu \cite{WX} classified symmetric graphs of prime order $p$,
of order $2p$ and of order $3p$, respectively.
For valency $3, 4$ and $5$ case, see \cite{CLP15,FZL16,GZF11,LWX,YFKL19,ZF10}
and references therein for examples.
Recently, a classification of 7-valent symmetric graphs of
order $2pq$ with $p,q$ distinct primes was given by
Hua, Chen and Xiang \cite{HLX},
and a characterization of 7-valent symmetric graphs
of order $4p^n$
was obtained in \cite{PY18}.

A typical method for studying
symmetric graphs is taking normal quotient graphs,
stated as following.
Let $\Ga$ be a $G$-arc-transitive graph.
For an intransitive normal subgroup
$N$ of $G$,
denote by $V\Ga_N$
the set of $N$-orbits in $V\Ga$. The {\it normal quotient graph}
$\Ga_N$ of $\Ga$ induced by $N$ is defined  with
vertex set $V\Ga_N$ and two vertices $B,C\in V\Ga_N$ are adjacent if and
only if some vertex in $B$ is adjacent in $\Ga$ to some vertex in
$C$. If $\Ga$ and $\Ga_N$ have the same valency, then
$\Ga$ is called a {\it normal cover} of $\Ga_N$.
In particular, $\Ga$ is called {\it basic}
if $\Aut\Ga$ has no nontrivial normal subgroup $N$ such that $\Ga$ and
$\Ga_N$ have the same valency.
From the definition, there is a natural `two-steps strategy' for studying symmetric graphs:
\vskip0.1in

Step 1. Determine their normal quotient graphs.

Step 2. Reconstruct the original graphs by
determining the normal covers of their normal quotient graphs.
\vskip0.1in

The main purpose of this paper is to determine
all the basic normal quotient graphs
of the connected $7$-valent symmetric graphs of order $2pq^n$
with $p<q$ odd primes and $n\ge 2$.

The notation used in this paper are standard,
see \cite{Atlas}.
For example, for a positive integer $n$, we use $\ZZ_n$ and $\D_{2n}$
to denote the cyclic group of order $n$ and the dihedral group of order $2n$,
respectively.
For two groups $N$ and $H$, denote by $N\times H$ the direct product of $N$ and $H$,
by $N.H$ an extension of $N$ by $H$, and
if such an extension is split, then we write $N:H$ instead of $N.H$.

For a positive integer $s$,
an $s$-{\it arc} of $\Ga$ is a sequence $v_0,v_1,\dots,v_s$ of $s+1$ vertices of $\Ga$
such that $v_{i-1},v_i$ are adjacent for $1\le i\le s$ and $v_{i-1}\ne v_{i+1}$ for $1\le i\le s-1$.
If $G\le\Aut\Ga$ is transitive on the set of $s$-arcs of $\Ga$,
then $\Ga$ is called {\it $(G,s)$-arc-transitive};
if $\Ga$ is $(G,s)$-arc-transitive
but not $(G,s+1)$-arc-transitive,
then $\Ga$ is called {\it $(G,s)$-transitive}.
Particularly, a $(\Aut\Ga,s)$-transitive graph is simply called
{\it $s$-transitive}.

Our main result is as follows. For convenience, graphs
appearing in Table 1 are introduced in Section 2.

\begin{theorem}\label{Thm-1}
Let $\Ga$ be a connected symmetric graph of valency $7$ and order $2pq^n$,
with $p<q$ odd primes and $n\ge 2$.
Then $\Ga$ is $s$-transitive with $1\le s\le 3$,
and is a normal cover of one of the basic graphs $\Sig$ listed in Table $1$.
\end{theorem}

\begin{table}
\[\begin{array}{llll|llll} \hline
\Sig & (p,q) & \Aut\Sig & s & \Sig & (p,q) & \Aut\Sig & s \\ \hline
\K_{7,7} & p=7 & \S_7\wr\ZZ_2 & 3 & \CC_{30} & (3,5) & \S_8 & 2  \\
\HS(50) & (3,5)& \PSU(3,5).\ZZ_2 & 2 & \CC_{78}^{1} & (3,13) & \PSL(2,13) & 1\\
\CC_{78}^{2} & (3,13) & \PGL(2,13) & 1 
& \CC_{310} & (5,31) & \Aut(\PSL(5,2)) & 3 \\
\CD(2p,7) & 7\mid p-1 & \D_{2p}\times\ZZ_7 & 1\\
\hline
\end{array}\]
\caption{Basic normal quotient graphs of symmetric graphs of order $2pq^n$}
\end{table}

For any given positive integer $k$,
a result of Conder, Li and Poto\v cnik \cite{CLP15}
asserts that
there are only finitely many connected $2$-arc-transitive $7$-valent graphs of
order $kp$ or $kp^2$ with $p$ a prime. Theorem~\ref{Thm-1}, together with \cite[Theorem 1.1]{PY18}
(for case $p=2$)
have the following corollary.

\begin{corollary}\label{Cor-1}
\begin{itemize}
\item[(1)] For any given positive integer $n$, there are only finitely many
connected $2$-arc-transitive $7$-valent graphs of order $2pq^n$ with $7\ne p<q$ primes.
\item[(2)] For any given positive integer $n$, there is no
connected $2$-arc-transitive $7$-valent graphs of order $2pq^n$ with $7<p<q$ primes.
\end{itemize}
\end{corollary}

\section{Preliminaries}

\subsection{Examples} 
As usual, for a positive integer $n$, denote by $\K_n,\K_{n,n}$ and $\K_{n,n}-n\K_2$
the complete graph of order $n$, the complete bipartite graph
of order $2n$, and the graph deleted a 1-match from $\K_{n,n}$, respectively.
Also, the Hoffman-Singleton graph of order $50$ and valency 7 is denoted by
$\HS(50)$.

For a group $G$ and a subset $S\subseteq G\setminus\{1\}$,
with $S=S^{-1}:=\{g^{-1}\mid
g\in S\}$, the {\it Cayley graph} of the group $G$ with respect to
$S$ is with vertex set $G$
and two vertices $g$ and $h$ are adjacent if and
only if $hg^{-1}\in S$. This Cayley graph is denoted by $\Cay(G,S)$.

\begin{example}\label{Dihedrant}
Let $G=\l a,b\mid a^m=b^2=1,a^b=a^{-1}\r\cong\D_{2m}$ be a dihedral group,
with $m$ a positive integer. Let $k$
be a solution of the congruence equation
$$ x^{6}+x^{5}+\cdots + x+1\equiv 0~(\mod~m).$$
Define a Cayley graph
$$\CD(2m,7)=\Cay(G,\{ b,ab,a^{k+1}b,\dots,a^{k^{5}+k^{4}+\cdots+1}b\}).$$
Then $\CD(2m,7)$ is a connected arc-transitive graph of valency $7$. In particular, if
$m\ge 31$, then $\CD(2m,7)$ is arc-regular
and $\Aut(\CD(2m,7))=\D_{2m}:\ZZ_7$,
see \cite[Theorem B; Proposition 4.1]{Feng-Li}.
\end{example}

The following are several specific examples, see \cite[Section 3]{HLX}.

\begin{example}\label{Exam-2}
\begin{itemize}
\item[(i)] There is unique connected symmetric $7$-valent graph of order $30$,
denoted by $\CC_{30}$, and $\Aut(\CC_{30})=\S_8$;

\item[(ii)] There are exactly two connected symmetric $7$-valent graphs of order $78$,
denoted by $\CC_{78}^1$ and $\CC_{78}^2$,
and $\Aut(\CC_{78}^1)=\PSL(2,23)$ and $\Aut(\CC_{78}^2)=\PGL(2,23)$;

\item[(iii)] There is unique connected symmetric $7$-valent graph of order $310$,
denoted by $\CC_{310}$, and $\Aut(\CC_{310})=\Aut(\PSL(5,2))$.

\end{itemize}
\end{example}

\begin{lemma}\label{graphs}
Let $\Ga$ be a connected $7$-valent symmetric graph.
Then the following statements hold, where $p<q$ are primes.
\begin{itemize}
\item[(1)] {\rm(\cite[Table]{C-O})} If $|V\Ga|=2p$, then
           $\Ga=\K_{7,7}$ or $\CD(2p,7)$ with $7\mid p-1$.
\item[(2)] {\rm(\cite[Section 4]{HLX})} If $|V\Ga|=2pq$, then
$\Ga=\CC_{30},\CC_{78}^{1},\CC_{78}^2,\CC_{310}$ or $\CD(2pq,7)$ with
$p=7$ or $7\mid p-1$ and $7\mid q-1$.
\end{itemize}
\end{lemma}

\subsection{Background results}

For a positive integer $m$ and a group $T$,
denote by $\pi(m)$ the number of the primes
which divide $m$, and by $\pi(T)$ the set of the primes
dividing $|T|$.
The group $T$ is called a $K_n$-group if $|\pi(T)|=n$.
The simple $K_n$-groups with $3\le n\le 6$ are classified in \cite{BW} and \cite{AA}.

\begin{theorem}\label{K_5-group}{\rm(\cite[Theorem A]{AA})}
Let $T$ be a simple $K_5$-group. Then one of the following holds:
\begin{itemize}
\item[(a)] $T=\PSL(2,q)$ with $\pi(q^2-1)=4$;
\item[(b)] $T=\PSU(3,q)$ with $\pi((q^2-1)(q^3+1))=4$;
\item[(c)] $T=\PSL(3,q)$ with $\pi((q^2-1)(q^3-1))=4$;
\item[(d)] $T=O_5(q)$ with $\pi(q^4-1)=4$;
\item[(e)] $T=\Sz(2^{2m+1})$ with $\pi((2^{2m+1}-1)(2^{4m+2}+1))=4$;
\item[(f)] $T=R(3^{2m+1})$ with $\pi((3^{4m+2}-1))=3$ and $\pi(3^{4m+2}-3^{2m+1}+1)=1$;
\item[(g)] $T=\A_{11}, \A_{12}, \M_{22},J_3,HS,He,McL,\PSL(4,4),\PSL(4,5),\PSL(4,7),
\PSL(5,2),\\
\PSL(5,3),\PSL(6,2),O_7(3), O_9(2), \PSp(6,3),
\PSp(8,2), \PSU(4,4),\PSU(4,5),\\
\PSU(4,7),\PSU(4,9),\PSU(5,3),\PSU(6,2),O^+(8,3),
O^{-}(8,2), ^3 D_4(3), G_2(4), G_2(5),\\
G_2(7)$ or $G_2(9)$.
\end{itemize}
\end{theorem}

The vertex stabilizers of connected 7-valent symmetric graphs were determined independently by
{\rm \cite[Theorem 1.1]{GSZ}} and {\rm \cite[Theorem 3.4]{LLW16}},
where $\F_n$ with $n$ a positive integer denotes the Frobenius group of order $n$.

\begin{lemma}\label{stabilizer}
Let $\Ga$ be a connected $7$-valent $(G,s)$-transitive graph, where $G\le\Aut\Ga$ and $s\ge 1$.
Then $s\le3$ and one of the following holds, where $\a\in V\Ga$.
\begin{itemize}
\item[(a)] If $G_{\a}$ is soluble, then $|G_{\a}|\mid 2^2\cdot3^2\cdot7$.
Further, the couple $(s,G_{\a})$ is listed in the following table.

\[\begin{array}{l|l|l|l} \hline
s & 1 & 2 & 3 \\ \hline
G_{\a} &  \ZZ_7,~\F_{14},~\F_{21},~\F_{14}\times\ZZ_2,~\F_{21}\times\ZZ_3& \F_{42},
~\F_{42}\times\ZZ_2,~\F_{42}\times\ZZ_3& \F_{42}\times\ZZ_6  \\ \hline
\end{array}\]
\vskip0.2in

\item[(b)] If $G_{\a}$ is insoluble, then $|G_{\a}|\mid 2^{24}\cdot 3^4\cdot 5^2\cdot 7$.
Further, the couple $(s,G_{\a})$ is listed in the following table.

\[\begin{array}{l|l|l} \hline
s & 2 & 3  \\ \hline
G_{\a} & \PSL(3,2), \ASL(3,2),&  \PSL(3,2){\times}\S_4,\A_7{\times}\A_6,\S_7{\times}\S_6,(\A_7{\times}\A_6){:}\ZZ_2,   \\
       & \ASL(3,2){\times}\ZZ_2,\A_7,\S_7 &
       \ZZ_2^6{:}(\SL(2,2){\times}\SL(3,2)),[2^{20}]{:}(\SL(2,2){\times}\SL(3,2))  \\ \hline
       |G_{\a}|&2^3{\cdot}3{\cdot}7,2^6{\cdot}3{\cdot}7,&2^6{\cdot}3^2{\cdot}7,~
       2^6{\cdot}3^4{\cdot}5^2{\cdot}7,~2^8{\cdot}3^4{\cdot}5^2{\cdot}7,~2^7{\cdot}3^4{\cdot}5^2{\cdot}7,\\
         & 2^7{\cdot}3{\cdot}7,2^3{\cdot}3^2{\cdot}5{\cdot}7,
         2^4{\cdot}3^2{\cdot}5{\cdot}7& 2^{10}{\cdot}3^2{\cdot}7,~~2^{24}{\cdot}3^2{\cdot}7   \\ \hline
\end{array}\]
\end{itemize}

In particular, if $5\mid |G_{\a}|$, then $|G_{\a}|\mid 2^{8}\cdot3^4\cdot 5^2\cdot 7$,
and $G_{\a}^{\Ga(\a)}\cong\A_7$ or $\S_7$;
if $5{\not |}~|G_{\a}|$, then $|G_{\a}|\mid 2^{24}\cdot3^2\cdot 7$.
\end{lemma}

The following theorem is a special case of
\cite[Lemma 2.5]{Li-Pan}
which slightly improves a nice result of Praeger \cite[Theorem 4.1]{Praeger92}.

\begin{theorem}\label{Praeger} Let $\Ga$ be a connected $G$-arc-transitive
graph of odd prime valency, and let $N\lhd G$ have more than two orbits on
$V\Ga$, where $G\le\Aut\Ga$. Then the following statements hold.
\begin{itemize}
\item[(1)] $N$ is semiregular on $V\Ga$, $G/N\le\Aut\Ga_N$,
$\Ga_N$ is $G/N$-arc-transitive, and $\Ga$ is a normal $N$-cover
of $\Ga_N$;
\item[(2)] $\Ga$ is $(G,s)$-arc-transitive
if and only if $\Ga_N$ is $(G/N,s)$-arc-transitive,
where $1\le s\le 5$ or $s=7$;
\item[(3)] $G_{\a}\cong(G/N)_{\d}$, where $\a\in V\Ga$ and $\d\in V\Ga_N$.
\end{itemize}
\end{theorem}

A transitive permutation group $X\le\Sym(\Ome)$ is called {\it quasiprimitive} if each
minimal normal subgroup of $X$ is transitive on $\Ome$, while $X$ is
called {\it biquasiprimitive} if each of its minimal normal
subgroups has at most two orbits and there exists one which has exactly two orbits
on $\Ome$.

We have a next generalization of \cite[Lemma 5.1]{FZL16}.

\begin{lemma}\label{2p^n}
Let $\Ga$ be a connected $G$-arc-transitive $r$-valent graph of order
$2q^n$, where $G\le\Aut\Ga$, $n\ge 2$, and $r\ge 5$ and $q\ge 5$  are primes.
Then either $\Ga=\K_{7,7}$ or $\HS(50)$, 
or $G$ has a minimal normal elementary abelian $q$-subgroup.
\end{lemma}

\proof If $G$ is quasiprimitive or biquasiprimitive on $V\Ga$,
by \cite[Theorem 1.2]{PL18}, Lemma~\ref{2p^n} is true.
Suppose that $G$ is neither quasiprimitive nor biquasiprimitive on $V\Ga$.
Then $G$ has a minimal normal subgroup $N$ which has at least three orbits on $V\Ga$,
by Theorem~\ref{Praeger}, $N$ is semiregular on $V\Ga$
and hence $|N|\mid 2q^n$.
It follows that either $N=\ZZ_2$ or $N=\ZZ_q^d$ for some $d<n$.
For the former case, the normal quotient graph $\Ga_N$ is arc-transitive of
odd order $q^n$ and odd valency $r$,
a contradiction.
Therefore $N=\ZZ_q^d$, as required.\qed

\section{Technical lemmas}

The two lemmas regarding simple groups in this section are based on the classifications
of simple $K_n$-groups with $3\le n\le 6$, obtained in \cite{BW} and \cite{AA}.

\begin{lemma}\label{5-factors}
Let $r<s$ be odd primes,
and let $T$ be a nonabelian simple group such that
$|T|\mid2^{25}\cdot 3^2\cdot7rs^{l}$ and $7rs^l\mid|T|$
for $l\geq 1$.
Then one of the following  holds.

\begin{itemize}
\item[(i)] $|\pi(T)|=4$, and $T$ is isomorphic to one of the groups listed in Table~$\ref{4-factor}$.

\begin{table}
\[\begin{array}{ll|ll|ll} \hline
T & |T| & T & |T| & T & |T|  \\ \hline
\J_2 & 2^7\cdot3^3\cdot5^2\cdot7 & \A_7 & 2^3\cdot3^2\cdot5\cdot7 & \A_8 & 2^6\cdot3^2\cdot5\cdot7  \\
\PSL(2,13) & 2^2\cdot3\cdot7\cdot13 & \PSL(2,27) & 2^2\cdot3^3\cdot7\cdot13 & \PSL(2,97) & 2^5\cdot3\cdot7^2\cdot97  \\
\PSL(2,127) & 2^3\cdot3^2\cdot7\cdot127  & \PSL(3,4) & 2^6\cdot3^2\cdot5\cdot7 &
\PSL(3,8) & 2^9\cdot3^2\cdot7^2\cdot73 \\
\PSU(3,5) & 2^4\cdot3^2\cdot5^3\cdot7  &&&& \\ \hline
\end{array}\]
\caption{Some simple $K_4$-groups}\label{4-factor}
\end{table}

\item[(ii)] $|\pi(T)|=5$, and $T$ is isomorphic to one of the groups listed in Table~$\ref{5-factor}$.
In particular, $l=1$.

\end{itemize}
\end{lemma}

\begin{table}
\[\begin{array}{l|ll} \hline
T & \M_{22} & \PSL(5,2)  \\
|T| & 2^7\cdot3^2\cdot5\cdot7\cdot11 & 2^{10}\cdot3^2\cdot5\cdot7\cdot31 \\ \hline

T & \M_{22} & \PSL(5,2)  \\
|T| & 2^7\cdot3^2\cdot5\cdot7\cdot11 & 2^{10}\cdot3^2\cdot5\cdot7\cdot31 \\ \hline

T & \PSL(2,2^6) & \PSL(2,29)  \\
|T| & 2^6\cdot3^2\cdot5\cdot7\cdot13 & 2^2\cdot3\cdot5\cdot7\cdot29 \\
\hline
T & \PSL(2,41) & \PSL(2,43)  \\
|T| & 2^3\cdot3\cdot5\cdot7\cdot41 & 2^{10}\cdot3^2\cdot5\cdot7\cdot31 \\
\hline
T & \PSL(2,71) & \PSL(2,83)  \\
|T| & 2^3\cdot3^2\cdot5\cdot7\cdot71 & 2^2\cdot3\cdot7\cdot41\cdot83 \\
\hline
T &\PSL(2,113) & \PSL(2,167)  \\
|T| & 2^4\cdot3\cdot7\cdot19\cdot113 & 2^3\cdot3\cdot7\cdot83\cdot167 \\
\hline
T & \PSL(2,223) & \PSL(2,503)  \\
|T| & 2^5\cdot3\cdot7\cdot37\cdot223 & 2^3\cdot3^2\cdot7\cdot251\cdot503 \\
\hline
T & \PSL(2,673) & \PSL(2,2017)  \\
|T| & 2^5\cdot3\cdot7\cdot337\cdot673 & 2^5\cdot3^2\cdot7\cdot1009\cdot2017 \\
\hline
T & \PSL(2,3583) & \PSL(2,64513)  \\
|T| &2^9\cdot3^2\cdot7\cdot199\cdot3583 & 2^{10}\cdot3^2\cdot7\cdot32257\cdot64513 \\
\hline
T & \PSL(2,2752513) & \PSL(2,16515073) \\
|T| & 2^{17}\cdot3\cdot7\cdot1376257\cdot2752513 & 2^{18}\cdot3^2\cdot7\cdot8257537\cdot16515073 \\

\hline
\end{array}\]
\caption{Some simple $K_5$-groups}\label{5-factor}
\end{table}

\proof Clearly, $3\le |\pi(T)|\le 5$.
If $|\pi(T)|=3$, by \cite[Theorem I]{BW},
there are exactly 8 specific simple $K_3$-groups listed in \cite[Table 1]{BW},
checking the orders, no group $T$ exists in this case.

\vskip0.1in
\noindent{(i).} Suppose $|\pi(T)|=4$. By \cite[Theorem I]{BW},
either
\vskip0.05in
\begin{itemize}
\item[(a)] $T$ is isomorphic to one of the groups listed in \cite[Table 2]{BW}; or
\item[(b)] $T=\PSL(2,q)$ for some prime power $q$.
\end{itemize}
\vskip0.05in

Assume (a) occurs.
Suppose $5\in \pi(T)$. As $7\in\pi(T)$,
by checking \cite[Table 2]{BW}, $T$ is a $\{2,3,5,7\}$-group,
and so $r,s\in \{3,5,7\}$.
If $s=7$, then $r=3$ or 5,
and one easily checks that
no group $T$ exists in the case.
If $s=5$, then $r=3$,
so $|T| \mid 2^{25}\cdot3^3\cdot5^l\cdot7$ and $7\cdot3\cdot5^l\mid |T|$,
and one may derive that
$T=\J_2,\A_7,\A_8,\PSL(3,4)$ and $\PSU(3,5)$.
Suppose now $5\notin \pi(T)$. By \cite[Table 2]{BW}, $r=3$ or 7 and $s>7$,
hence $|T| \mid 2^{25}\cdot3^3\cdot7^2\cdot s^l$,
by checking the orders,
we obtain $T=\PSL(3,8)$.

Now assume (b) occurs.
If $q$ is a power of $2,3$ or $7$,
by \cite[TABLE 3]{BW}, the only example is $T=\PSL(2,27)$.
For the other cases, by \cite[Theorem 3.2]{BW},
$q\ge 11$ is a prime, 
notice that $|\PSL(2,q)|$ is always divisible by $3$,
we conclude that $T$ is a $\{2,3,7,q\}$-group,
hence $s=q$, $l=1$, and $r\in\{3,7\}$.
Further, since $|T|\mid2^{25}\cdot 3^2\cdot7rs^{l}$,
we have ${q-1\over 2}\cdot{q+1\over 2}\mid2^{24}\cdot3^3\cdot7^2$,
and as $({q-1\over 2},{q+1\over 2})=1$,
it follows that either $q+1\mid2\cdot3^3\cdot7^2$ if
$2\mid {q-1\over 2}$,
or $q-1\mid2\cdot3^3\cdot7^2$ if $2\mid {q+1\over 2}$.
Then a computation by Magma~\cite{Magma} shows $q\in\{13,17,19,41,43,53,97,127,293,379,881,883\}$.
Checking the orders, we obtain $T=\PSL(2,13),\PSL(2,97)$ or $\PSL(2,127)$.

\vskip0.1in
\noindent{(ii).} Suppose $|\pi(T)|=5$. Then $T$ is a $\{2,3,7,r,s\}$-group
and satisfies part (a)-(g) of Theorem~\ref{K_5-group}.
We analyse these cases one by one in the following.
Notice that $|T| \mid 2^{25}\cdot3^2\cdot7rs^l$,
we obtain
\begin{eqnarray}\label{eq1}
2^{26}{\not |}~|T|,~~3^3{\not |}~|T|,~~7^2{\not |}~|T|,~~r^2{\not |}~|T|.
\end{eqnarray}

Assume part (a) of Theorem~\ref{K_5-group} occurs.
Then
\begin{eqnarray}\label{eq2}
|T|=|\PSL(2,q)|={1\over (2,q-1)}q(q-1)(q+1).
\end{eqnarray}
Since $|T|\mid 2^{25}\cdot3^2\cdot7rs^l$
and $7rs^l\mid |T|$,
we have $q\neq 3,3^2,7$ or $r$,
and then derive from Eq.(\ref{eq1})
that $q=2^i$ with $1\leq i \leq 25$ or $q=s^l$.
For the former case, since $\pi(q^2-1)=4$,
we get $q=2^6,2^8,2^9,2^{11}$ or $2^{23}$,
and by checking the orders, we obtain $T=\PSL(2,2^6)$.
For the latter case, we have ${q+1\over 2}\cdot {q-1\over 2}\mid2^{24}\cdot3^2\cdot7r$,
and as $({q+1\over 2},{q-1\over 2})=1$,
it follows that either $q-1\mid2^{25}\cdot3^2\cdot7$ if $r\mid{q+1\over 2}$,
or $q+1\mid2^{25}\cdot3^2\cdot7$ if $r\mid{q-1\over 2}$.
Recall that $\pi(q^2-1)=4$, $q=s^l$ with $s>7$ and
$T$ satisfies Eq.(\ref{eq1}), a direct computation by Magma\cite{Magma}
shows that
$q=29,41,43,71,83,113,167,223,503,673,2017,3583,64513,2752513$ or $16515073$,
as in Table 3.

Assume part (b) occurs.
Then
\begin{eqnarray}\label{eq3}
|T|=|\PSU(3,q)|={1\over (3,q+1)}q^3(q-1)(q+1)^2(q^2-q+1).
\end{eqnarray}
By Eq.(\ref{eq1}), $q$ is a 2-power or a $s$-power.
If $q$ is a 2-power, then $q=2^i$ with $1\leq i \leq 8$,
and as $\pi((q^2-1)(q^3+1))=4$,
we get $q=2^4,2^5$ or $2^7$;
however, in these three cases,
$|T|$ is always not divisible by $7$, a contradiction.
If $q$ is an $s$-power, as $7^2{\not |}~|T|$ and $r^2{\not |}~|T|$,
we obtain $(q+1)^2\mid2^{25}\cdot 3^2$,
or equivalently $q+1\mid2^{12}\cdot3$.
Since $\pi((q^2-1)(q^3+1))=4$, computation in Magma \cite{Magma} shows $q=11$ and 23;
however, in both cases, $|T|$ is not divisible by $7$,
also a contradiction.

Assume part (c) occurs. Then
\begin{eqnarray}\label{eq4}
|T|=|\PSL(3,q)|={1\over (3,q-1)}q^3(q-1)^2(q+1)(q^2+q+1).
\end{eqnarray}
By Eq.(\ref{eq1}), one derives $q$ is a 2-power or a $s$-power.
Then with similar discussion as in part (b) above,
one may draw a contradiction.

Assume part (d) occurs.
Then
\begin{eqnarray}\label{eq5}
|T|=|O_5(q)|={1\over 2}q^4(q^4-1)(q^3-1)(q^2-1).
\end{eqnarray}
By Eq.(\ref{eq1}), $q$ is a 2-power or a $s$-power.
If $q$ is a 2-power, then $q=2^i$ with $1\leq i \leq 6$.
Since $\pi(q^4-1)=4$, we have $q=2^3$ or $2^4$,
and $|T|=2^{12}\cdot3^4\cdot5\cdot7^2\cdot13$ or $2^{16}\cdot3^2\cdot5^2\cdot17^2\cdot257$
respectively,
contradicting Eq.(\ref{eq1}).
If $q$ is a $s$-power, as $7^2{\not |}~|T|$ and $r^2{\not |}~|T|$,
we have $(q^2-1)^2\mid2^{25}\cdot3^2$, and so ${q+1\over 2}\cdot{q-1\over 2}\mid2^{11}\cdot3$.
Noting that $({q+1\over 2},{q-1\over 2})=1$,
we conclude that either ${q+1\over 2}\mid3$ or ${q-1\over 2}\mid3$,
implying $s\le q\le 7$, a contradiction.

Assume part (e) occurs.
Then $|T|=|\Sz(2^{2m+1})|=2^{4m+2}(2^{4m+2}+1)(2^{2m+1}-1)$.
Since $2^{26}{\not |}~|T|$ and $\pi((2^{2m+1}-1)(2^{4m+2}+1))=4$,
we derive $m=3$; however, $|\Sz(2^7)|=2^{14}\cdot5\cdot29\cdot113\cdot127$
is not divisible by 7, a contradiction.

Assume part (f) occurs.
Then $|T|=|R(3^{2m+1})|=3^{6m+3}(3^{6m+3}+1)(3^{2m+1}-1)$,
so $3^9\mid |T|$, contradicting $3^3{\not |}~|T|$.

Finally, assume part (g) occurs.
Checking the orders of the 30 specific simple groups there,
we obtain
$T=\M_{22}$ and $\PSL(5,2)$.\qed

\begin{lemma}\label{6-factors}
Let $r<s$ be odd primes,
and let $T$ be a nonabelian simple group such that $|T| \mid 2^9\cdot 3^4\cdot5^2\cdot7rs^{l}$
and $35rs^l\mid |T|$
with $l\geq 1$.
Then one of the following holds.

\begin{itemize}
\item[(i)] $|\pi(T)|=4$, and $T$ is isomorphic to one of the groups listed in Table~$\ref{4-factor1}$.

\begin{table}
\begin{center}
\begin{tabular}{llllll} \hline
  $T$ & $|T|$ & $T$ & $|T|$ & $T$ & $|T|$   \\ \hline
  $\J_2$ & $2^7\cdot3^3\cdot5^2\cdot7$ & $\A_{10}$ & $2^7\cdot3^4\cdot5^2\cdot7$ & $\PSU(3,5)$ & $2^4\cdot3^2\cdot5^3\cdot7$  \\
  $\PSp(4,7)$ & $2^8\cdot3^2\cdot5^2\cdot7^4$ & $\PSL(2,49)$ & $2^4\cdot3\cdot5^2\cdot7^2$ & & \\ \hline
\end{tabular}
\end{center}
\caption{Certain simple $K_4$-groups} \label{4-factor1}
\end{table}

\item[(ii)] $|\pi(T)|=5$, $l=1$ and $T$ is isomorphic to one of the groups listed in Table~$\ref{5-factor1}$.

\begin{table}
\begin{center}
\begin{tabular}{llllll} \hline
  $T$ & $|T|$ & $T$ & $|T|$   \\ \hline
  $\A_{11}$ & $2^7\cdot3^4\cdot5^2\cdot7\cdot11$ & $\A_{12}$ & $2^9\cdot3^5\cdot5^2\cdot7\cdot11$   \\
  $\M_{22}$ & $2^7\cdot3^2\cdot5\cdot7\cdot11$ & $\HS$ & $2^9\cdot3^2\cdot5^3\cdot7\cdot11$ \\
  $\PSL(2,2^6)$ & $2^6\cdot3^2\cdot5\cdot7\cdot13$ & $\PSL(2,5^3)$ & $2^3\cdot3^2\cdot5^3\cdot7\cdot31$ \\
  $\PSL(2,29)$ & $2^2\cdot3\cdot5\cdot7\cdot29$ & $\PSL(2,41)$ & $2^3\cdot3\cdot5\cdot7\cdot41$  \\
  $\PSL(2,71)$ & $2^3\cdot3^2\cdot5\cdot7 \cdot71$ & $\PSL(2,251)$ & $2^2\cdot3^2\cdot5^3\cdot7\cdot251$ \\
  $\PSL(2,449)$ & $2^6\cdot3^2\cdot5^2\cdot7\cdot449$ &   &   \\ \hline
\end{tabular}
\end{center}
\caption{Certain simple $K_5$-groups} \label{5-factor1}
\end{table}

\item[(iii)] $|\pi(T)|=6$, $l=1$ and $T=\J_1,\M_{23}$, or $\PSL(2,q)$
with $q=139,181,211,239,281,\\
349,379,421,601,631,701,769,811,
839,1009,1049,1051,1399,1511,1889,2099,\\
2239,2267,2269,2591,2689,2801,3779,4481,6481,
6719,7559,10079,12601,15121,\\21601,26881,28351,30241,37799,53759,56701,
69119,96769,172801,201599,\\
453599,483839$ or $907199$.
\end{itemize}
\end{lemma}

\proof Obviously, $3\le |\pi(T)|\le 6$.
If $|\pi(T)|=3$, by {\rm \cite[Theorem I]{BW}},
$T$ is isomorphic to one of the 8 groups listed in {\rm \cite[Table 1]{BW}}.
However, the order of each group there is not divisible by $35$,
a contradiction.

\vskip0.1in
\noindent{(i).} Assume $|\pi(T)|=4$.
By \cite[Theorem I]{BW},
$T$ is isomorphic to one of the groups listed in \cite[Table 3]{BW} or
$T=\PSL(2,q)$ for some prime power $q$.
For the former case,
since $35rs^l\mid |T|$
and $|T| \mid 2^9\cdot 3^4\cdot5^2\cdot7rs^{l}$,
one easily derives that
$T=\J_2,\A_{10},\PSU(3,5)$ or $\PSp(4,7)$.
For the later case, notice that $3\mid |\PSL(2,q)|$ and $35\mid |T|$,
$T$ is a $\{2,3,5,7\}$-group,
then by \cite[TABLE 3]{BW}, the only example is
$T=\PSL(2,49)$.

\vskip0.1in
\noindent{(ii).}  Assume $|\pi(T)|=5$.
Then $s>7$ and $T$ satisfies parts (a)-(g) of Theorem~\ref{K_5-group}.
Since $|T|\mid 2^9\cdot 3^4\cdot5^2\cdot7rs^{l}$, we have
\begin{eqnarray}\label{eq6}
2^{10}{\not |}~|T|,~~3^6{\not |}~|T|,~~5^4{\not |}~|T|,~~7^3{\not |}~|T|.
\end{eqnarray}
Suppose $T=\PSL(2,q)$, as in part (a) of Theorem~\ref{K_5-group}.
Then $T$ is a $\{2,3,5,7,s\}$-group as $3\mid |\PSL(2,q)|$.
If $q$ is a 2-power, then $q=2^6,2^8$ or $2^9$ since $2^{10}{\not |}~|T|$ and $\pi(q^2-1)=4$,
by checking the orders, we obtain $T=\PSL(2,2^6)$.
If $q$ is a 3-power, then $q\in \{3,3^2,3^3,3^4,3^5\}$ since $3^6{\not |}~|T|$,
it follows that $\pi(q^2-1)\ne 4$, a contradiction.
If $q$ is a 5-power, then $q=5^3$ because $5^4{\not |}~|T|$ and $\pi(q^2-1)=4$,
which gives rise to an example $T=\PSL(2,5^3)$.
If $q$ is a 7-power, then $q=7$ or $7^2$ since $7^3 {\not |}~|T|$,
contradicting $\pi(q^2-1)=4$.
Now, assume that $q$ is a $s$-power.
Then ${q+1\over 2}\cdot{q-1\over 2}\mid2^9\cdot 3^5\cdot5^3\cdot7^2$.
Since $({q+1\over 2},{q-1\over 2})=1$, we have
${q-1\over 2}\mid3^5\cdot5^3\cdot7^2$ or ${q+1\over 2}\mid3^5\cdot5^3\cdot7^2$.
Recall that $\pi(q^2-1)=4$, computation in Magma\cite{Magma} shows
$q\in \{29,41,43,71,89,149,151,251,269,271,293,449,751,809,2251,
2647,4051,7937,12149,\\
20249,23813\}$.
Checking the orders, we obtain the examples $T=\PSL(2,29),
\PSL(2,41)$ and $\PSL(2,449)$.

Suppose $T=\PSU(3,q)$, as in part (b).
Since $\pi((q^2-1)(q^3+1))=4$, by Eq.(\ref{eq3}) and (\ref{eq6}),
we derive
that $q$ is a $s$-power, and $(q+1)^2\mid2^{10}\cdot3^5\cdot5^3\cdot7^2$,
so $q+1\mid2^5\cdot3^2\cdot5\cdot7$.
Since $\pi((q^2-1)(q^3+1))=4$, a computation by Magma \cite{Magma} shows
that $q\in \{11,13,17,19,23\}$;
however, by checking the orders, no group $T$ exists in the case.
Similarly, one may exclude part (c), namely $T=\PSL(3,q)$.

Suppose $T=\Sz(2^{2m+1})$ or $R(3^{2m+1})$,
as in part (d) or (e).
Then $|T|=2^{4m+2}(2^{4m+2}\\
+1)(2^{2m+1}-1)$ or $3^{6m+3}(3^{6m+3}+1)(3^{2m+1}-1)$,
respectively.
Since $|\pi(\Sz(2^3))|=4$, $T\ne\Sz(2^3)$,
and hence $2^{10}\mid|T|$ or $3^9\mid|T|$,
contradicting Eq.(\ref{eq6}).

Suppose $T=O_5(q)$, as in part (f).
Since $\pi(q^4-1)=4$, by Eq.(\ref{eq5}) and (\ref{eq6}), we conclude that
$q$ is a $s$-power, and $(q-1)^3\mid2^{10}\cdot3^5\cdot5^3$,
hence $q-1\mid2^3\cdot3\cdot5$.
It follows $q=11$ and $|T|=|O_5(11)|=2^8\cdot3^2\cdot5^2\cdot11^4\cdot61$
is not divisible by $7$, a contradiction.

Finally, suppose $T$ lies in the groups listed in part (g).
Checking the orders,
we obtain $T=\A_{11},\A_{12},\M_{22},$ or $\HS$.

\vskip0.1in
\noindent{(iii).} Assume that $|\pi(T)|=6$. Then $7<r<s$ and $s>11$.
By {\rm \cite[Theorem B]{AA}},
one of the following holds:

\begin{itemize}
\item[(a)] $T=\PSL(2,q)$ where $\pi(q^2-1)=5$;
\item[(b)] $T=\PSL(3,q)$ where $\pi((q^2-1)(q^3-1))=5$;
\item[(c)] $T=\PSL(4,q)$ where $\pi((q^2-1)(q^3-1)(q^4-1))=5$;
\item[(d)] $T=\PSU(3,q)$ where $\pi((q^2-1)(q^3+1))=5$;
\item[(e)] $T=\PSU(4,q)$ where $\pi((q^2-1)(q^3+1)(q^4-1))=5$;
\item[(f)] $T=O_5(q)$ where $\pi(q^4-1)=5$;
\item[(g)] $T=G_2(q)$ where $\pi(q^6-1)=5$;
\item[(h)] $T=\Sz(2^{2m+1})$ where $\pi((2^{2m+1}-1)(2^{4m+2}+1))=5$;
\item[(i)] $T=R(3^{2m+1})$ where $\pi((3^{2m+1}-1)(3^{6m+3}+1))=5$;
\item[(j)] $T$ is one of the 38 groups listed in \cite[Theorem B]{AA}.
\end{itemize}

Recall that $|T|\mid 2^9\cdot 3^4\cdot5^2\cdot7rs^{l}$,
then we have
\begin{eqnarray}\label{eq7}
2^{10}{\not |}~|T|,3^5{\not |}~|T|,5^3{\not |}~|T|,7^2{\not |}~|T|,r^2{\not |}~|T|.
\end{eqnarray}

Since
$|\Sz(2^{2m+1}|=2^{4m+2}(2^{4m+2}+1)(2^{2m+1}-1)$
and $2^{10}{\not |}~|T|$, we have $m=1$
and $|\pi(T)|=|\pi(\Sz(8))|=4\ne 6$,
this contradiction excludes case (h).
Since $|R(3^{2m+1})|=3^{6m+3}(3^{6m+3}+1)(3^{2m+1}-1)$,
$3^9\mid |T|$, contradicting Eq.(\ref{eq7}),
this excludes case (i).
For case (j),
by checking the orders,
we have $T=\J_1$ or $\M_{23}$.

Suppose case (a) occurs.
Since $\pi(q^2-1)=5$, $q\ne 3,3^2,3^3,3^4,5,5^2,7$,
and by Eq.(\ref{eq7}), we have that either $q=2^i~(1\leq i \leq 9)$ or $s^l$.
The former case does not give examples by checking the orders.
For the latter case, by Eq.(\ref{eq2}),
we have ${q-1\over 2}\cdot{q+1\over 2}\mid2^8\cdot 3^4\cdot5^2\cdot7r$,
and as $({q-1\over 2},{q+1\over 2})=1$, it follows
that either ${q+1\over 2}\mid2^9\cdot 3^4\cdot5^2\cdot7$
or ${q-1\over 2}\mid2^9\cdot 3^4\cdot5^2\cdot7$.
Recall that $\pi(q^2-1)=5$ and $s>11$, computation in Magma \cite{Magma}
shows that $q$ lies in part (iii) of Lemma~\ref{6-factors}.

Suppose case (b) occurs.
Since $\pi((q^2-1)(q^3-1))=5$,
$q\ne 2,2^2,2^3$ or $3$,
and by Eq.(\ref{eq4}) and (\ref{eq7}),
we derive that $q$ is a $s$-power,
and $(q-1)^2\mid2^9\cdot 3^4\cdot5^2$,
implying $q-1\mid2^4\cdot 3^2\cdot5$.
Since $\pi((q^2-1)(q^3-1))=5$, a computation by Magma \cite{Magma}
shows that $q=37,41$ or 241,
which does not give rise to examples by checking the orders.
Similarly, one may exclude case (d).

Suppose case (e) occurs.
Then
$$|T|=|\PSU(4,q)|={1\over (4,q+1)}q^6(q^2-1)^2(q^2+1)(q^3+1).$$
By Eq.(\ref{eq7}), one may conclude that $q$ is a $s$-power,
and $(q^2-1)^2\mid2^9\cdot 3^4\cdot5^2$,
hence ${q-1\over 2}\cdot{q+1\over 2}\mid2^2\cdot 3^2\cdot5$.
As $({q-1\over 2},{q+1\over 2})=1$,
we obtain
${q-1\over 2}\mid 2^2\cdot5$ or ${q+1\over 2}\mid 2^2\cdot5$.
It follows that $q=19$ since $\pi((q^2-1)(q^3+1)(q^4-1))=5$,
and $|T|=|\PSU(4,19)|=2^7\cdot3^4\cdot5^3\cdot7^3\cdot19^6\cdot181$,
contradiction Eq.(\ref{eq7}).
For case (c), then
$$|T|=|\PSL(4,q)|={1\over (4,q-1)}q^6(q^2-1)(q^3-1)(q^4-1),$$
a similar argument as in case (e) may draw a contradiction.

Suppose case (f) occurs.
Since $\pi(q^4-1)=5$, $q\ne 2,2^2$ and $3$,
then by Eq.(\ref{eq5}) and (\ref{eq7}), we conclude that $q$ is a $s$-power,
and $(q^2+1)(q^3-1)(q^2-1)^2\mid2^{10}\cdot 3^4\cdot5^2\cdot7r$,
hence $(q^2-1)^2\mid2^{10}\cdot 3^4\cdot5^2$,
or equivalently $q^2-1\mid2^5\cdot 3^2\cdot5$.
As discussed in the above paragraph,
one may derive that
$q=17,19,29,31$ or $89$,
which does not give rise to example by checking the orders.
Finally, for case (g), $|T|=|G_2(q)|=q^6(q^6-1)(q^2-1)$,
a similar discussion may draw a contradiction.\qed

\section{Vertex quasiprimitive and vertex biquasiprimitive cases}

Let $\Ga$ be a connected $G$-arc-transitive 7-valent graph of order $2pq^n$,
where $G\le\Aut\Ga$, $p<q$ are odd primes and $n\geq2$.
Let $N$ be a minimal normal subgroup of $G$. Then
$N=T^d$, with $T$ a simple group and $d\ge 1$.
Let $\a \in V\Ga$.

\begin{lemma}\label{d<2}
If $N$ is nonabelian, then $d=1$.
\end{lemma}
\proof Suppose for a contradiction that $N$ is nonabelian
and $d\ge 2$.
Then $|N|{\not |}~2pq^n$, $N_{\a}\ne 1$,
and $N$ has at most two orbits on $V\Ga$ by Theorem~\ref{Praeger}.
Set $N=T_1\times T_2\times\cdots\times T_d$
with each $T_i\cong T$.

Assume first $N$ is transitive on $V\Ga$.
Since $1\ne N_{\a}\lhd G_{\a}$
and $\Ga$ is connected,
we have
$1\ne N_{\a}^{\Ga(\a)}\lhd G_{\a}^{\Ga(\a)}$.
It follows that $N_{\a}^{\Ga(\a)}$ is transitive,
and $\Ga$ is $N$-arc-transitive.
If $T_1$ is transitive on $V\Ga$,
then the centralizer $C_N(T_1)$ is semiregular on $V\Ga$ (see \cite[Theorem 4.2A]{DM}),
so is $T_2$, which is a contradiction as $|T_2|{\not |}$ does not divide $|V\Ga|=2pq^n$;
if $T_1$ has at least three orbits on $V\Ga$,
by Theorem~\ref{Praeger}, $T_1$ is semiregular,
again a contradiction.
Therefore, $T_1$ has exactly 2 orbits, say $U$ and $W$, on $V\Ga$.
Since $T_1\lhd N$, $U$ and $W$ form a $N$-block system on $V\Ga$.
It follows that the set stabilizer $N_U$ is of index 2 in $N$,
which is a contradiction because $N=T^d$ has no subgroup with index 2.

Assume now $N$ has exactly two orbits, say $\Del_1$ and $\Del_2$, on $V\Ga$.
Then $\Ga$ is a bipartite graph with bipartitions $\Del_1$ and $\Del_2$.
Let $G^+=G_{\Del_1}=G_{\Del_2}$, the stabilizer on the bipartitions.
If $G^+$ acts unfaithfully on $\Del_1$,
by \cite[Lemma 5.2]{GLP03}, $\Ga$ is a complete bipartite graph,
so $\Ga=\K_{7,7}$ as $\val(\Ga)=7$ and hence $|V\Ga|=14$,
a contradiction.
Suppose $G^+$ acts faithfully on $\Del_1$.
Then $N\le G^+$ can be viewed as a transitive permutation group on $\Del_1$.
If $T_1$ is transitive on $\Del_1$,
then \cite[Theorem 4.2A]{DM} implies $T_2$ is semiregular on $\Del_1$,
hence $|T_2|\mid pq^n$, a contradiction.
Thus $T_1$ has at least two orbits on $\Del_1$. It then follows from
\cite[Lemma 3.2]{LWX} that $T_1$ is semiregular on $\Del_1$,
also a contradiction.\qed

The next two lemmas exclude the vertex quasiprimitive and vertex biquasiprimitive cases.

\begin{lemma}\label{quasiprimitive}
If $G$ is quasiprimitive on $V\Ga$, then no graph $\Ga$ exists.
\end{lemma}

\proof Since $G$ is quasiprimitive on $V\Ga$,
$N$ is transitive on $V\Ga$.
If $N$ is abelian, then
$N$ is regular on $V\Ga$
and so $|T|^d=|N|=2pq^n$,
a contradiction.
Thus $N$ is nonabelian, and then by
Lemma~\ref{d<2},
we have $d=1$ and $N=T$.
Further, since $T_{\a}\ne 1$,
we conclude that $\Ga$ is $T$-arc-transitive,
and hence $T_{\a}$ satisfies Lemma~\ref{stabilizer}.
We divided our proof into two cases depending on whether 5 divides $|T_{\a}|$ or not.

\vskip0.1in
{\noindent\bf Case 1.} Assume $5{\not |}~|T_\a|$.

By Lemma~\ref{stabilizer}, $|T_{\a}|\mid 2^{24}\cdot3^2\cdot 7$,
and by the transitivity of $T$,
we have $|T|=|V\Ga||T_{\a}|$ divides $2^{25}\cdot3^2\cdot7pq^n$;
on the other hand, since $\Ga$ is $T$-arc-transitive, we have $7\mid |T_{\a}|$,
and so $7pq^n\mid |T|$.
Therefore, $T$ satisfies Lemma~\ref{5-factors};
in particular, $|\pi(T)|=4$ or 5.
If $|\pi(T)|=5$, by Lemma~\ref{5-factors}(ii),
$n=1$, a contradiction.

Suppose $|\pi(T)|=4$. Noting that $n\ge 2$, by Lemma~\ref{5-factors}(i),
we easily conclude that the couple $(T,p,q^n)=(\J_2,3,5^2)$ or
$(\PSU(3,5),3,5^3)$.
For the former case, $|V\Ga|=2\cdot3\cdot5^2$,
so $|T_\a|={|T| \over |V\Ga|}=2^6\cdot3^2\cdot7$;
however, by Atlas \cite{Atlas}, $\J_2$ has no subgroup with order $2^6\cdot3^2\cdot7$,
a contradiction.
For the latter case, $|V\Ga|=2\cdot3\cdot5^3$,
hence $|T_\a|={|T|\over |V\Ga|}=2^3\cdot3\cdot7$,
then by Lemma~\ref{stabilizer}, we obtain $T_\a=\PSL(3,2)$;
however, computation in Magma \cite{Magma} shows that
no graph $\Ga$ exists in this case.

\vskip0.1in
{\noindent\bf Case 2.} Assume $5\mid|T_\a|$.

By Lemma~\ref{stabilizer},
$|T_{\a}|\mid 2^{8}\cdot3^4\cdot 5^2\cdot 7$,
and $T_{\a}^{\Ga(\a)}\cong\A_7$ or $\S_7$.
It follows that
$|T|=|V\Ga||T_{\a}|$ divides $2^9\cdot3^4\cdot5^2\cdot7pq^n$;
moreover, as $35\mid|T_\a|$, we have $35pq^n$ divides $|T|$.
Therefore, $T$ satisfies Lemma~\ref{6-factors}; in particular $4\le\pi(T)\le 6$.
If $|\pi(T)|=5$ or 6, by Lemma~\ref{6-factors},
we have $n=1$, a
contradiction.

Suppose $|\pi(T)|=4$. Since $n\ge 2$, by Lemma~\ref{6-factors}(i),
the only possibility is $T=\PSp(4,7)$ and $(p,q)=(3,7^3)$ or $(5,7^3)$.
Consequently, $|V\Ga|=2\cdot3\cdot7^3$ or $2\cdot5\cdot7^3$,
and $|T_\a|={|T| \over |V\Ga|}=2^7\cdot3\cdot5^2\cdot7$ or $2^7\cdot3^2\cdot5\cdot7$,
respectively.
By Lemma~\ref{stabilizer}, it is a contradiction.\qed

\begin{lemma}\label{biquasiprimitive}
If $G$ is biquasiprimitive on $V\Ga$, then no graph $\Ga$ exists.
\end{lemma}

\proof Since $G$ is biquasiprimitive on $V\Ga$,
$G$ has a minimal normal subgroup $N=T^d$ which has exactly two orbits
(say $\Del_1$ and $\Del_2$) on $V\Ga$.
Then $\Ga$ is a bipartite graph with bipartition $\Del_1$ and $\Del_2$.
Let $G^+=G_{\Del_1}=G_{\Del_2}$.
Then $N\le G^+$, $|G:G^+|=2$ and $G_{\a}=G^+_{\a}$.
If $N$ is abelian, then $N$ is regular on $\Del_1$
and so $|T|^d=|N|=pq^n$, a contradiction.
Hence $N$ is nonabelian, and by Lemma~\ref{d<2},
we further conclude that
$N=T$ is a nonabelian simple group.

If $G^+$ acts unfaithfully on $\Del_1$ or $\Del_2$,
by \cite[Lemma 5.2]{GLP03}, $\Ga$ is a complete bipartite graph,
so $\Ga=\K_{7,7}$ and $|V\Ga|=14$,
a contradiction.

Assume now $G^+$ acts faithfully on $\Del_1$ and $\Del_2$.
Then by \cite[Theorem 1.5]{LPAZ}), either
\begin{itemize}
\item[(1)] $G^+$ is quasiprimitive on $\Del_i$; or
\item[(2)] $G^+$ has two normal subgroups $M_1$ and $M_2$ such that $M_1\cong M_2$
           are semiregular on $V\Ga$. Further, the group $M_1\times M_2$ is regular on $\Del_i$.
\end{itemize}

For case (2), we have $|M_1|^2=|\Del_i|=pq^n$, a contradiction.

Suppose case (1) occurs.
Since $G^+$ is quasiprimitive on $\Del_i$ and has a
simple minimal normal subgroup $T$,
by O'Nan-Scott-Praeger theorem (\cite{Praeger92}),
$\soc(G^+)=T$ or $T^2$.
For the latter case, $G^+$ is of holomorph type
and $T$ is regular on $\Del_i$,
so $|T|=pq^n$, a contradiction.
Therefore, $\soc(G^+)=T$.
Further,
if $T$ is not the unique minimal normal subgroup of
$G$, since $G=G^+.\ZZ_2$,
one easily derives
$G=G^+\times\ZZ_2$,
hence the normal subgroup $\ZZ_2$ has $pq^n$ orbits on $V\Ga$,
contradicting the biquasiprimitivity of $G$.
Thus $G$ is almost simple with socle $T$,
and we may set $G=T.o$, and $G^+=T.{o'}$ with
$\ZZ_2\le o\le\Out(T)$ and $|o:o'|=2$.

\vskip0.1in
{\noindent\bf Case 1.} Assume $5{\not |}~|T_\a|$.

Since $T_{\a}\le G_{\a}$,
by Lemma~\ref{stabilizer},
$|T_{\a}|\mid 2^{24}\cdot3^2\cdot 7$,
and hence $|T|=|\Del_1||T_{\a}|$ divides $2^{24}\cdot3^2\cdot7pq^n$;
on the other hand, noting that $T_{\a}\ne 1$, we obtain $7\mid |T_{\a}|$,
and so $7pq^n\mid |T|$.
Therefore, $T$ satisfies Lemma~\ref{5-factors},
and $\pi(T)=4$ or $5$.

If $\pi(T)=5$, by Lemma~\ref{5-factors}(ii), we have $n=1$, a contradiction.

Suppose now $|\pi(T)|=4$. Then
by Lemma~\ref{5-factors}(i) and notices that $n\ge 2$,
we have $(T,p,q^n)=(\J_2,3,5^2)$ or
$(\PSU(3,5),3,5^3)$.
For the former case,
as $\Out(\J_2)\cong\ZZ_2$,
we have $o=\ZZ_2$, $o'=1$
and $G^+=T$.
It follows $|G_\a|=|G^+_\a|=|T_\a|={|T| \over |\Del_1|}=2^7\cdot3^2\cdot7$,
which is a contradiction by Lemma~\ref{stabilizer}.
For the latter case,
$|T_\a|={|T| \over |\Del_1|}=2^4\cdot3\cdot5\cdot7$.
Since $\Out(\PSU(3,5))\cong\S_3$,
we have that either $o=\S_3$ and $o'=\ZZ_3$, or $o=\ZZ_2$ and $o'=1$.
Thus, $|G_\a|=|G^+_\a|=|T_\a|\cdot|o'|=2^4\cdot3\cdot5\cdot7$ or $2^4\cdot3^2\cdot5\cdot7$.
By Lemma~\ref{stabilizer}, the only possibility is
$|G_\a|=2^4\cdot3^2\cdot5\cdot7$,
and $G_\a\cong S_7$;
however, computation in Magma \cite{Magma}
shows that no graph $\Ga$ exist in the case.

\vskip0.1in
{\noindent\bf Case 2.} Assume $5\mid|T_\a|$.

As $T_{\a}\lhd G_{\a}$, by Lemma~\ref{stabilizer},
$|T_{\a}|\mid 2^{8}\cdot3^4\cdot 5^2\cdot 7$,
and $T_{\a}^{\Ga(\a)}\cong\A_7$ or $\S_7$.
It follows that
$|T|=|V\Ga||T_{\a}|$ divides $2^8\cdot3^4\cdot5^2\cdot7pq^n$;
moreover, as $35\mid|T_\a|$, we have $35pq^n$ divides $|T|$.
Hence $T$ satisfies Lemma~\ref{6-factors}.

If $|\pi(T)|=5$ or 6, by Lemma~\ref{6-factors},
$n=1$, a contradiction.

Suppose $|\pi(T)|=4$. By Lemma~\ref{6-factors}(i),
we have $T=\PSp(4,7)$,
and $(p,q^n)=(3,7^3)$ or $(5,7^3)$,
so $|T_\a|={|T| \over |\Del_1|}=2^8\cdot3\cdot5^2\cdot7$ or $2^8\cdot3^2\cdot5^2\cdot7$,
respectively.
Further, as $\Out(\PSp(4,7))\cong\ZZ_2$,
we have $o=\ZZ_2$, $o'=1$ and $G^+=T.o'=T$.
Thus $|G_\a|=|G^+_\a|=|T_\a|=2^8\cdot3\cdot5^2\cdot7$ or $2^8\cdot3^2\cdot5^2\cdot7$,
by Lemma~\ref{stabilizer}, which is a contradiction.\qed

\section{Proofs of Theorems~\ref{Thm-1}}

We will complete the proof of Theorem~\ref{Thm-1} in this final section.

\begin{lemma}\label{norm-q-sg}
Let $p<q$ be odd primes, and let $\Ga$ be a connected $G$-arc-transitive $7$-valent graph
of order $2pq^m$, where $G\leq \Aut\Ga$ and $m\geq 2$.
Then either $\Ga$ is a $\ZZ_3$-cover of $\HS(50)$,
or $G$ has a normal elementary abelian $q$-subgroup.
\end{lemma}

\proof By Lemmas~\ref{quasiprimitive} and ~\ref{biquasiprimitive},
$G$ is neither quasiprimitive nor biquasiprimitive on $V\Ga$.
Hence $G$ has a minimal normal subgroup, say $N$,
which has at least three orbits on $V\Ga$.
By Theorem~\ref{Praeger}, $N$ is semiregular,
and hence $|N|$ divides $|V\Ga|=2pq$.
It follows that $N$ is soluble
and $N\cong\ZZ_2,\ZZ_p$ or $\ZZ_q^s$ with $s\le n$.
For the last case, we are done.
For the first case, by Theorem~\ref{Praeger},
$\Ga_N$ is connected arc-transitive of
odd order $pq^n$ and odd valency $7$,
a contradiction.
For the second case,
again by Theorem~\ref{Praeger},
$\Ga_N$ is $G/N$-arc-transitive of
order $2q^m$ and valency 7.
Clearly, $\Ga_N\ne\K_{7,7}$.
It then follows from Lemma~\ref{2p^n}
that either $\Ga_N=\HS(50)$, or
$G/N$ has a minimal normal subgroup $M/N\cong\ZZ_q^k$
for some positive integer $k$.
For the former case, $(p,q)=(3,5)$,
and $\Ga$ is a $\ZZ_3$-cover of $\HS(50)$.
For the latter case,
$M=\ZZ_p.\ZZ_q^k$,
and as $p<q$,
the Sylow $q$-subgroup $\ZZ_q^k$ of $M$ is characteristic in $M$,
and hence normal in $G$.\qed

Now, we are ready to prove Theorems~\ref{Thm-1}
Corollary~\ref{Cor-1}.

\vskip0.1in
{\noindent\bf Proof of Theorems~\ref{Thm-1}.}
Suppose $\Ga$ is $G$-arc-transitive with $G\le\Aut\Ga$.
Let $M$ be a maximal normal $q$-subgroup of $G$.
Clearly, $M$ has at least $2p\ge 6$ orbits on $V\Ga$.
Then by Theorem~\ref{Thm-1},
$\Ga$ is a normal cover of $\Ga_M$,
and $\Ga_M$ is a connected $G/M$-arc-transitive
graph of order $2pq^{m}$ with $0\le m\le n-1$.

If $m=0$, then $\Ga_M$ is of order $2p$,
by Lemma~\ref{graphs}(1),
$\Ga=\K_{7,7}$ or $\CD_{2p}$ with $7\mid p-1$.

If $m=1$, then $\Ga_M$ is of order $2pq$,
by Lemma~\ref{graphs}(2),
$\Ga_M=\K_8,\CC_{30},\CC_{78}^{1},\CC_{78}^2,\CC_{310}$,
or $\CD(2pq,7)$ with $7\mid q-1$
For the last case,
as $2pq>31$, by Example~\ref{Dihedrant},
$\Ga_M=\CD(2pq,7)$ is arc-regular,
so $G/M=\Aut(\Ga_M)\cong\D_{2pq}:\ZZ_7$.
Now, $G$ has a normal subgroup
$H=M.\ZZ_q$ which has $2p$ orbits on $V\Ga$,
and $\Ga$ is a normal cover of $\Ga_H=\CD(2p,7)$ by Theorem~\ref{Thm-1}.

Now assume $m\ge 2$. By Lemma~\ref{norm-q-sg},
either $\Ga$ is a $\ZZ_3$-cover of $\HS(50)$,
or $G$ has a normal elementary abelian $q$-subgroup, say $X/M$.
For the former case, $\Ga$ is a normal $M.\ZZ_3$-cover of $\HS(50)$.
For the latter case, $X$ is a normal $q$-subgroup of $G$,
and by the maximality of $M$, $X$ has at most two orbits on $V\Ga$.
It follows that $pq^n\mid |X|$, a contradiction.

Finally, one easily verifies that the graphs in Table 1 are basic graphs.
This completes the proof of Theorems~\ref{Thm-1}.\qed

\vskip0.1in
{\noindent\bf Proof of Corollary~\ref{Cor-1}.}
Let $\Ga$ be a connected $2$-arc-transitive $7$-valent graphs of order $2pq^n$
with $q>p$ primes.

If $n=1$, notice that $\CD(2pq,7)$ is not 2-arc-transitive,
Corollary~\ref{Cor-1} is true by Lemma~\ref{graphs}(2).

Assume $n\ge 2$. If $p\ne 7$, by Theorems~\ref{Thm-1} and \ref{Praeger},
$\Ga$ is a normal cover of $\CC_{30}$, $\HS(50)$ or $\CC_{310}$,
hence $(p,q)=(3,5)$ or $(5,31)$.
Now Corollary~\ref{Cor-1} easily follows.\qed

\end{document}